\documentstyle[amssymb, amstex,12pt]{article}

\def\p{\partial}
\def\w{\wedge}

\def\dim{\operatorname{dim}}

\def\ct{\operatorname{ct}}
\def\ft{\operatorname{ft}}
\newcommand{\CC}{{\Bbb C}}

\newcommand{\ZZ}{{\Bbb Z}}
\newcommand{\QQ}{{\Bbb Q}}
\newcommand{\lan}{\langle}
\newcommand{\ran}{\rangle}
\renewcommand{\a}{\alpha}
\renewcommand{\b}{\beta}
\renewcommand{\c}{\gamma}
\renewcommand{\d}{\delta}
\newcommand{\gel}{\varepsilon}

\newcommand{\gl}{\lambda}

\newcommand{\go}{\omega}

\newcommand{\gt}{\tau}

\newcommand{\calh}{{\cal H}}

\newcommand{\calm}{{\cal M}}

\newcommand{\calo}{{\cal O}}

\def\ev{\operatorname{ev}}
\def\z{{\bf z}}

\title{On the WDVV-equation in quantum K-theory}
\author{Alexander Givental}
\date{UC Berkeley and Caltech}
\begin{document}
\maketitle

{\bf 0. Introduction.} 
Quantum cohomology theory can be described in general words as intersection
theory in spaces of holomorphic curves in a given K\"ahler or almost
K\"ahler manifold $X$. 
By quantum $K$-theory we may similarly understand a study of
complex vector bundles over the spaces of holomorphic curves in $X$. 
In these notes, we will introduce a $K$-theoretic
version of the Witten-Dijkgraaf-Verlinde-Verlinde equation which 
expresses the associativity constraint of the ``quantum multiplication"
operation on $K^*(X)$.

Intersection indices of cohomology theory,
\[ \int _{[\ \text{space of curves}\ ]} \go_1\w ... \w\go_k \]
obtained by evaluation on the fundamental cycle of cup-products of 
cohomology classes, are to be replaced in $K$-theory by Euler
characteristics 
\[ \chi (\ \text{space of curves}\ ; V_1\otimes ... \otimes V_k) \]
of tensor products of vector bundles.
The hypotheses needed in the definitions of the intersection indices and
Euler characteristics --- that the spaces of curves are compact and non-singular, 
or that the bundles are holomorphic --- are rarely satisfied.
We handle this foundational problem by restricting ourselves throughout 
the notes to the setting where the problem disappears. 
Namely, we will deal with the so called 
{\em moduli spaces $X_{n,d}$ of degree $d$ genus
$0$ stable maps to $X$ with $n$ marked points} 
{\bf assuming that $X$ is a homogeneous K\"ahler space}. 
Under the hypothesis, the moduli spaces $X_{n,d}$ (we will review 
their definition and properties when needed) are known
to be compact complex orbifolds (see \cite{Kn, BM}). We use
their fundamental cycle $[X_{n,d}]$, well-defined over $\QQ $, in the 
definition of intersection indices, and we use sheaf cohomology
in the definition of the Euler characteristic of a holomorphic {\em orbi-bundle} $V$:  
\[ \chi (X_{n,d}; V) := \sum (-1)^k \dim H^k (X_{n,d}; \Gamma (V)) .\]     

\medskip

{\bf 1. Correlators.} The WDVV-equation is usually formulated 
in terms of the following generating function for {\em correlators}: 
\[ F(t,Q)= \sum_{d} \sum_{n=0}^{\infty} \frac{Q^d}{n!} (t,...,t)_{n,d} .\]
Here $d\in H_2(X,\ZZ)$ runs the Mori cone of {\em degrees}, that is homology classes 
represented by fundamental cycles of rational holomorphic curves 
in $X$, and the correlators $(\phi_1,...,\phi_n)_{n,d}$ are defined using
the {\em evaluation maps} at the marked points:
\[ \ev_1\times ... \times \ev_n : X_{n,d} \to X\times ... \times X .\]
In cohomology theory, we pull-back to the moduli space $X_{n,d}$ the 
$n$ cohomology classes $\phi_1,...,\phi_n \in H^*(X,\QQ)$ of $X$
and define the correlator among them by 
\[ (\phi_1,...,\phi_n)_{n,d}:=\int_{[X_{n,d}]} \ev_1^*(\phi_1)\w ...\w
\ev_n^*(\phi_n) .\]
In K-theory, we pull-back $n$ elements $\phi_1,...,\phi_n\in K^*(X)$ 
(representable under our restriction on $X$ by holomorphic vector bundles 
or their formal differences) and put
\[ (\phi_1,...,\phi_n)_{n,d}:=\chi (X_{n,d} ; \ev_1^*(\phi_1)\otimes ...
\otimes \ev_n^*(\phi_n) ) .\] 
We will treat the series $F$ as a formal function of $t\in H$ 
depending on formal parameters $Q=(Q_1,...,Q_{\text{Betti}_2(X)})$, 
where $H=H^*(X,\QQ)$ or $H=K^*(X)$.   

Let $\{\phi_{\a}\}$ be a graded basis in $H^*(X,\QQ)$, and 
\[ g_{\a\b}:= \lan \phi_{\a},\phi_{\b}\ran=\int_{[X]} \phi_{\a}\w\phi_{\b}
\] denote the intersection matrix. Let $(g^{\a\b})=(g_{\a\b})^{-1}$ be the
inverse matrix
(so that $\sum (\phi_{\a}\otimes 1) g^{\a\b} (1\otimes \phi_{\b})$ is 
Poincare-dual to the diagonal in $X\times X$). In quantum cohomology theory, 
one defines the {\em quantum cup-product} $\bullet $ on the tangent space $T_tH$ by 
\[ \lan \phi_{\a}\bullet \phi_{\b},\phi_{\c} \ran := F_{\a\b\c}(t) \]
(where the subscripts on the RHS mean partial derivatives in the basis
$\{ \phi_{\a} \}$). In the above notation the associativity of
the quantum cup-product is equivalent to the WDVV-identity:

\medskip

{\em $\sum_{\gel,\gel'} F_{\a\b\gel}g^{\gel\gel'}F_{\gel'\c\d}$ 
is totally symmetric in $\a,\b,\c,\d$.}

\medskip

{\bf 2. Stable maps, gluing and contraction.} In order to explain the
proof of the WDVV-identity we have to discuss some properties of
the moduli spaces $X_{n,d}$ (see \cite{Kn, BM, FP} for more details). 

We consider prestable marked curves $(C, \z)$, that is compact connected complex
curves $C$ with at most double singular points and with $n$ marked points
$\z=(z_1,...,z_n)$ which are non-singular and distinct. Two holomorphic
maps, $f:(C,\z)\to X$ and $f':(C',\z')\to X$, are called {\em equivalent}
if they are identified by an isomorphism $(C,\z)\to (C',\z')$ of the
curves. This definition induces the concept of {\em automorphism} of
a map $f: (C,\z)\to X$, and one calls $f$ {\em stable} if it has no 
non-trivial infinitesimal automorphisms. The moduli spaces $X_{n,d}$
consist of equivalence classes of stable maps with fixed number $n$
of marked points, degree $d$ and arithmetical genus $0$ (it is defined 
as $g=\dim H^1(C, \calo_C)$). 

In plain words, the space of degree $d$ holomorphic spheres in $X$ with
$n$ marked points is compactified by prestable curves which are trees of 
$\CC P^1$'s and satisfy the stability condition: each irreducible component 
$\CC P^1$ mapped to a point in $X$ must carry at least $3$ marked or singular points. 
Under the hypothesis that $X$ is a homogeneous K\"ahler
space, the moduli space $X_{n,d}$ has the structure of a compact complex orbifold of 
dimension $\dim_{\CC}X + \int_{d} c_1(T_X) + n -3$.  

In the case when $X$ is a point the moduli spaces coincide with the Deligne-Mumford
compactifications $\bar{\calm }_{0,n}$ of moduli spaces of configurations of marked
points on $\CC P^1$. For instance, $\calm_{0,4}$ is the set $\CC P^1-\{ 0,1,\infty \} $
of legitimate values of the cross-ratio of 4 marked points on $\CC P^1$.  
The compactification $\bar{\calm}_{0,4} =\CC P^1$ fills-in the fibidden values of
the cross-ratio by equivalence classes of the reducible curves $\CC P^1\cup \CC P^1$ 
with one double point and two marked point on each irreducible component. 
 
For $n\geq 3$, there is a natural {\em contraction} map $X_{n,d} \to \bar{\calm}_{0,n}$
defined by composing the map $f: (C,\z)\to X$ with $X\to pt$ (so that the components
of $C$ carrying $<3$ special points become unstable) and contracting the unstable 
components. Similarly, one can define the {\em forgetting} maps $\ft_i: X_{n+1,d}\to X_{n,d}$
by disregarding the $i$-th marked point and contracting the component if it has become
unstable. 

In particular, we will make use of the contraction map 
\[ \ct:  X_{n+4,d}\to \bar{\calm}_{0,4} \]
defined by forgetting the map $f: (C,\z) \to X$ and all the marked points except the first four.       
A legitimate value $\gl=\ct [f]$ of the cross-ratio means the following: 
the curve $C$ has a component $C_0=\CC P^1$ carrying $4$ special points with the cross-ratio $\gl$,
and the first $4$ marked point are situated on the branches of the tree connected to
$C_0$ at those $4$ special points. A forbidden value $\ct [f] =0$,$1$ or $\infty $ means that
$C$ containing a {\em chain} $C_0,...,C_k$ of $k>0$ of $\CC P^1$'s such that $2$ of the $4$ branches
of the tree carrying the marked points are connected to the chain via $C_0$, and the other two 
--- via $C_k$. Such stable maps form a stratum of codimension $k$ in the moduli space $X_{n,d}$.
We will refer to them as strata (or stable maps) {\em of depth $k$}.

A stable map of depth $1$ is glued from $2$ stable maps obtained by disconnecting 
$C_0$ from $C_1$. This gives rise to the {\em gluing map} 
\[ X_{n_0+3,d_0}\times_{\Delta} X_{n_1+3,d_1} \to X_{n_0+n_1+4,d_0+d_1}\] 
as follows. Consider the map from $X_{n_0+3,d_0}\times X_{n_1+3,d_1}$
to $X\times X$ defined by evaluation at the $3$-rd marked points. The source of
the gluing map is the preimage of the diagonal $\Delta \subset X\times X$. 
\footnote{Note that for a homogeneous K\"ahler $X$, the evaluation map is 
conveniently transverse to the diagonal in $X\times X$.} 
It consists of pairs of stable maps which have the same image of the third marked point 
and which therefore can be glued at this point into a single stable map of degree
$d_0+d_1$ with $n_0+2+n_1+2$ marked points.   

Similarly, gluing stable maps of depth $k$ from $k+1$ stable maps subject to
$k$ diagonal constraints at the double points of the chain $C_0,...,C_k$
defines appropriate gluing maps parameterizing the strata of depth $k$.

\medskip   

{\bf 3. Proof of the WDVV-identity.} All points in $\bar{\calm}_{0,4}$
represent the same (co)homology class. Thus the analytic fundamental cycles
of the fibers $\ct^{-1}(\gl )$ are homologous in $X_{n+4,d}$. The cohomological
WDVV-identity follows from the fact that for $\gl=0$, $1$ or $\infty $
the fiber $\ct^{-1} (\gl)$ consists of strata of depth $>0$, and moreover ---
the corresponding gluing maps (for all splittings $d=d_0+d_1$ of the degree
and all splittings of the $n=n_0+n_1$ marked points), being isomorphisms at
generic points, identify the analytic fundamental cycle of the fiber with
the sum of the fundamental cycles of $X_{n_0+3,d_1}\times_{\Delta }X_{n_1+3,d_2}$.   
This allows one to equate $3$ quadratic expression of the correlators
which differ by the order of the indices $\a,\b,\c,\d$ associated with the 
first $4$ marked points. 

We leave the reader to work out some standard 
combinatorial details which are needed in order to translate this argument into
the WDVV-identity for the generating function $F$ and note only that the contraction
with the intersection tensor $(g^{\gel\gel'})$ in the WDVV-equation takes care 
of the diagonal constraint $\Delta\subset X\times X$ for the evaluation maps.    

In K-theory, similarly, the push-forward to $X\times X$ of the structural sheaf
$\calo_{\Delta}$ of the diagonal is expressed as
\[ \sum (\phi_{\gel }\otimes 1)g^{\gel\gel'}(1 \otimes \phi_{\gel'}) \]  
via $(g^{\gel,\gel'})$ inverse to the ``intersection matrix"
\[ g_{\a\b}:=\lan \phi_{\a},\phi_{\b}\ran = \chi (X; \phi_{\a}\otimes \phi_{\b}) .\]
The argument justifying the WDVV-equation fails, however, since the above gluing map 
to $\ct^{-1}(\gl )$ is one-to-one only at the points of depth $1$ and does not 
identify the corresponding structural sheaves. Indeed, a stable map of depth $k$ 
can be glued from two stable maps in $k$ different ways and thus belongs to the
$k$-fold self-intersection in the image of the gluing map.

Let us examine the variety $\ct^{-1}(\gl)$ at a point of depth $k>1$. One of the
properties of Kontsevich's compactifications $X_{m,d}$ is that {\em after passing 
to the local non-singular covers} (defined by the orbifold structure of the
moduli spaces) {\em the compactifying strata form a divisor with normal crossings} 
\cite{Kn, BM}. Moreover, analyzing (inductively in $k$) the local structure of
the contraction map $\ct: X_{n+4,d} \to \bar{\calm}_{0,4}$ near a depth-$k$ point,
one easily finds the local model $\gl (x_1,...,x_k,...) = x_1...x_k$ for the
map $\ct $ in a suitable local coordinate system. In this model, the components 
$x_1=0,...,x_k=0$ of the divisor with normal crossings represent the strata of depth $1$, 
their intersections $x_{i_1}=x_{i_2}=0$ --- the strata of depth $2$, etc.
Denote by $\calo$ the algebra of functions on our local chart, so that  
$\calo / (x_{i_1},...,x_{i_l})$, $i_1<...<i_l$, are the algebras of functions on the depth-$l$ strata.
We have the following exact sequence of $\calo$-modules:  
\[ 0 \to \calo / (x_1...x_k) \to  \oplus\calo /(x_i) \to \oplus \calo /(x_{i_1},x_{i_2})
\to \oplus \calo /(x_{i_1},x_{i_2},x_{i_3}) \to ... .\]
Notice that the $\oplus$-terms in the sequence are the algebras of functions on the normalized
strata of depth $1$, depth $2$, etc. 
Translating this local formula to a global K-theoretic statement about gluing maps, we
conclude that in the Grothendieck group of orbi-sheaves on $X_{n+4,d})$, 
the element represented by the structural 
sheaf of $\ct^{-1}(\gl)$ for $\gl =0$,$1$ or $\infty$ is identified with the structural sheaf 
of the corresponding alternated disjoint sum over positive depth strata: 
\[ \sum X_{n_0+3,d_0}\times_{\Delta} X_{n_1+3,d_1} - 
\sum X_{n_0+3,d_0}\times_{\Delta} X_{n_1+2,d_1}\times_{\Delta}X_{n_2+3,d_2}+...\]

\medskip

{\bf 4. Formulation and consequences.}
Now we can apply the above K-theoretic statement about the moduli spaces 
to our generating functions.

Introduce 
\[ G(t,Q):=\frac{1}{2} \sum_{\a,\b} g_{\a\b}t_{\a}t_{\b}\ +\ F(t,Q) .\]
Let $(G^{\a\b})$ be the matrix inverse to $(G_{\a\b})=(\p_{\a}\p_{\b}G)$.

\medskip

{\bf Theorem.} 
\[ \sum_{\gel,\gel'} G_{\a\b\gel}G^{\gel\gel'}G_{\gel\c\d} \ 
\text{\em is totally symmetric in $\a,\b,\c,\d$.} \]   

\medskip

{\em Proof.} We have rewritten 
\[ F_{\a\b\gel}g^{\gel\gel'}F_{\gel'\c\d}-
F_{\a\b\gel}g^{\gel\mu}F_{\mu\mu'}g^{\mu'\gel'}F_{\gel'\c\d}+...\]
using the famous matrix identity $1-F+F^2-...=(1+F)^{-1}$. $\square $
  
\medskip

Introduce the {\em quantum tensor product} on $T_tH$ (with $H=K^*(X)$) by
\[ ( \phi_{\a}\bullet \phi_{\b},\phi_{\c}):=G_{\a\b\c} (t), \]  
and the metric $( , )$ on $TH$ is defined by $(\phi_{\mu},\phi_{\nu}):= G_{\mu\nu}(t) $.

\medskip

{\bf Corollary 1.} {\em The operations $( , )$ and $\bullet$ define on the tangent 
bundle the structure of a formal commutative associative Frobenius algebra with the unity
$1$.} \footnote{At $t=0, Q=0$ it turns into the usual multiplicative structure on $K^*(X)$.}

\medskip

{\em Proof.} As in the cohomology theory, it is a formal corollary of the Theorem, except
that the statement about the unity $1$ means that $G_{\a,1,\b}=G_{\a\b}$ and follows 
from the simplest instance of the {\em string equation} in the K-theory: 
$(1, t,...,t)_{n+1,d}=(t,...,t)_{n,d}$. The last equality is obvious. Indeed,
the push-forward of the constant sheaf $1$ along the map $\ft: X_{n+1,d}\to X_{n,d}$
forgetting the first marked point is the constant sheaf $1$ on $X_{n,d}$ since
the fibers are curves $C$ of zero arithmetic genus, $g=\dim H^1(C,\calo_{C})=0$, while
$H^0(C,\calo_{C})=\CC $ by Liouville's theorem. $\square $

\medskip   

We introduce on $T^*H$ the $1$-parametric family of connection operators
\[ \nabla_{q}:= (1-q) d - \sum_{\a} (\phi_{\a} \bullet )\ dt_{\a}\w\ .\]

{\bf Corollary 2.} {\em The connections $\nabla_{q}$ are flat for any $q\neq 1$.} 

\medskip

{\em Proof.} This follows from $\phi_{\a}\bullet \phi_{\b}=\phi_{\b}\bullet \phi_{\a}$, 
$d^2=0$, and $\p_{\a} (\phi_{\b}\bullet)=\p_{\b}(\phi_{\a}\bullet)$:
\[ \p_{\a} (\phi_{\b}\bullet )_{\mu}^{\nu}= 
G_{\mu\a\b\gel}G^{\gel\nu}-G_{\mu\b\gel}G^{\gel\gel'}G_{\gel'\a\gel''}G^{\gel''\nu} \]
is symmetric with respect to $\a$ and $\b$ due to the WDVV-identity. $\square $

\medskip

{\bf Proposition.} {\em The operator $\nabla_{-1}$ is twice the Levi-Civita connection
of the metric $(G^{\a\b})$ on $T^*H$.}

\medskip

{\em Proof.} For a metric of the form $G_{\a\b}=\p_{\a}\p_{\b}G$ the famous explicit 
formulas for the Christoffel symbols yield  
\[ 2\Gamma_{\a\b}^{\c}=[ G_{\a\gel,\b}+G_{\b\gel,\a}-G_{\a\b,\gel}]G^{\gel\c}=
G_{\a\b\gel}G^{\gel\c} = (\phi_{\b}\bullet)_{\a}^{\c}.\] 

\medskip

{\bf Corollary 3.} {\em The metric $( , )$ on $TH$ is flat.}

\medskip

We complete this section with a description of flat sections of the
connection operator $\nabla_q$ in terms of K-theoretic ``gravitational descendents". 
Let us introduce the generating functions
\[ S_{\a\b}(t,Q):=g_{\a\b}+\sum_{n,d} \frac{Q^d}{n!}
 ( \phi_{\a},t,...,t,\frac{\phi_{\b}}{1-qL} )_{n+2,d} ,\]
where the correlators are defined by
\[ (\psi_1,...,\psi_nL^k)_{m,d}:=\chi (X_{m,d};
\ev_1^*(\psi_1)\otimes ... \otimes \ev_m^*(\psi_m)\otimes
L^{\otimes k}) .\]
Here $L$ is the line {\em orbi}bundle over the moduli space 
$X_{m,d}$ of stable maps $(C,\z)\to X$ formed by the cotangent lines to $C$ at the {\em last} marked point 
(as specified by the position of the geometrical series
$1+qL+q^2L^2+...=(1-qL)^{-1}$ in the correlator).

\medskip

{\bf Theorem.} {\em The matrix $S:=(S_{\mu\nu})$ is a fundamental solution
to the linear PDE system:}
\[ (1-q)\p_{\a} S = (\phi_{\a}\bullet ) S \ . \]

\medskip

{\em Proof.} Taking $\phi_{\mu},\phi_{\a},\phi_{\b}$ and 
$\phi_{\nu}/(1-qL)$ for the content of the four distinguished marked
points in the proof of the WDVV-identity, we obtain its generalization
in the form:
\[ G_{\mu\a\gel}G^{\gel\gel'}\p_{\b}S_{\gel'\nu} = 
G_{\mu\b\gel}G^{\gel\gel'}\p_{\a}S_{\gel'\nu} ,\]
or $(\phi_{\a}\bullet ) \p_{\b}S = (\phi_{\b}\bullet ) \p_{\a} S$. 
Now it remains to put $\phi_{\b}=1$ and use $(1-q)\p_1 S =S$, which is
another instance of the string equation: 
\[ (1,t,...,t,\phi L^k)_{n+2,d} = 
(t,...,t,\phi (1+L+...+L^k))_{n+1,d}.\] 
The last relation is obtained by computing the push forward
of $L^{\otimes k}$ along $\ft_{1}:X_{n+2,d}\to X_{n+1,d}$.
\footnote{Some details can be found in \cite{W,P,L,G}. 
Briefly, one identifies the fibers of $\ft_1$ with the curves underlying the stable maps $f:(C,\z)\to X$ 
with $n+1$ marked points.
It is important to realize that the pull-back $L':=\ft_1^*(L)$ of the line
bundle named $L$ on $X_{n+1,d}$ differs from the line bundle named $L$ 
on $X_{n+2,d}$. In fact, there is a holomorphic section of 
$Hom (L',L)$ with the divisor $D$ defined by the last marked point 
$z_{n+1} \in C$, and the bundle $L$ restricted to $D$ is trivial 
(while $L'|_{D}$ is therefore conormal to $D$). Since $L'$ is trivial
along the fibers $C$, we find that $H^1(C,L^k)=0$ and
$H^0(C,L^k)=(L')^k \otimes H^0(C,\calo_{C}(kD))\simeq 
(L')^k(1+(L')^{-1}+...+(L')^{-k})$.   }

\medskip

{\bf 5. Some open questions.} 

{\em (a) Definitions.} It is natural to expect that the above
results extend from the case of homogeneous K\"ahler spaces $X$ to 
general compact K\"ahler and, even more generally, 
almost K\"ahler target manifolds.

In the K\"ahler case, the moduli of stable
degree $d$ genus $g$ maps with $n$ marked points form compact complex 
orbi-spaces $X_{g,n,d}$ equipped with the {\em intrinsic normal
cone} \cite{LT}. The cone gives rise \cite{F} to an element in K-group of 
$X_{g,n,d}$ which should be used in the definition of 
$K$-theoretic correlators in the same manner as the virtual fundamental cycle 
$[X_{g,n,d}]$ is used in quantum cohomology theory.

The moduli space $X_{g,n,d}$ can be also described as the zero locus
of a section of a bundle $E\to B$ over a non-singular space.
Due to the famous ``deformation to the normal cone" \cite{F}, 
the virtual fundamental cycle represents the Euler class of the bundle. 
This description survives in the almost K\"ahler case and yields a
topological definition and symplectic invariance of the
cohomological correlators. In $K$-theory, there exists a topological
construction of the push forward from $B$ to the point based on Whitney 
embedding theorem and Thom isomorphisms. We don't know however how to 
adjust the construction to our actual setting where $B$ is non-singular 
only in the {\em orbi}fold sense. 

One (somewhat awkward) option is
to define $K$-theoretic correlators topologically by the RHS of the
Kawasaki-Riemann-Roch-Hirzebruch formula \cite{Kw} for orbi-bundles
over $B$. This proposal deserves further study even in the K\"ahler
case since it may lead to a "quantum Riemann-Roch formula".             

\medskip

{\em (b) Frobenius-like structures.} Our results in Section $4$ show
that $K$-theoretic Gromov-Witten invariants of genus $0$ define 
on the space $H=K^*(X)$ a geometrical structure very similar to the
{\em Frobenius structure} \cite{D} of cohomology theory, but not identical
to it. 

One of the lessons is that the metric tensor on $H$, 
which can be in both cases described as $F_{\a,1,\b}$, is constant
in cohomology theory and equal to $g_{\a\b}$ only by an "accident", but
remains flat in $K$-theory even though it is not constant anymore. 

The translation $t\mapsto t+\gt 1$ in the direction of $1\in H$ 
leaves the structure invariant in cohomology theory, but causes 
multiplication by $e^{\gt}$ in $K$-theory --- because of a new form 
of the string equation.
Also, the $\ZZ$-grading missing in $K$-theory makes an important difference.
It would be interesting to study the axiomatic structure that emerges 
here and to compare it with the structure implicitly encoded by
$K$-theory on Deligne-Mumford spaces. 
 
\medskip

{\em (c) Deligne-Mumford spaces.} When the target space $X$ is the point,
the moduli spaces $X_{g,n,0}$ are Deligne-Mumford compactifications 
of the moduli spaces of genus $g$ Riemann surfaces with $n$ marked points.
The parallel between cohomology and $K$-theory suggest several 
problems.

Holomorphic Euler characteristics of universal cotangent line bundles
and their tensor products satisfy the string and dilation equations.
\footnote{The same is true not only for $X=pt$ (see \cite{Le}). By the
way, the push forward $\ft_*(L)$ along $\ft :X_{g,n+1,d}\to X_{g,n,d}$,
described by the dilation equation, equals $\calh+\calh^*-2+n$.
Here $\calh$ is the $g$-dimensional {\em Hodge bundle} with the fiber 
$H^1(C,\calo_{C})$. This answer replaces a similar factor $2g-2+n$ in the 
cohomological dilation equation, but also shows that tensor powers of
$\calh$ must be included to close up the list of ``observables".} 
$K$-theoretic generalization of the rest of Witten -- Kontsevich 
intersection theory \cite{W,K} is unclear. 

The case of genus $0$ and $1$ has been studied in \cite{P,L}
and \cite{Le}. The formula  
\[ \chi (\bar{\calm}_{0,n}; \frac{1}{(1-q_1L_1)...(1-q_nL_n)}) =
\frac{(1+q_1/(1-q_1)+...+q_n/(1-q_n))^{n-3}}
{(1-q_1)...(1-q_n)} \]
found by Y.-P. Lee \cite{L}
is analogous to the famous intersection theory result \cite{W, Kn}
\[ \int_{[\bar{\calm}_{0,n}]} \frac{1}{(1-x_1c_1(L_1))...(1-x_nc_1(L_n))} 
=(x_1+...+x_n)^{n-3} .\]
The latter formula is a basis for fixed point computations \cite{Kn, G}
in equivariant cohomology of the moduli spaces $X_{n,d}$ for toric
$X$. As it was notices by Y.-P. Lee, the former formula is not sufficient
for similar fixed point computation in K-theory: 
it requires Euler characteristics accountable for {\em invariants
with respect to permutations of the marked points}. 
Finding an $S_n$-equivariant version of Lee's formula is an important
open problem.  
   
\medskip

{\em (d) Computations.} The quantum $K$-ring is unknown even for
$X=\CC P^1$. It turns out that the WDVV-equation is not powerful enough
in the absence of grading constraints and {\em divisor equation} 
(see, for instance, \cite{G}). 

On the other hand, for $X=\CC P^n$, it is not hard  
to compute the generating functions $G(t,Q)$ and even $S_{\a\b}(t,Q,q)$ 
{\em at $t=0$} (see \cite{Le}). In cohomology theory, this would 
determine the {\em small} quantum cohomology ring due to the divisor
equation which, roughly speaking, identifies the $Q$-deformation 
at $t=0$ with the $t$-deformation at $Q=1$ along the subspace 
$H^2(X,\QQ) \subset H$. No replacement for the divisor equation 
seems to be possible in K-theory.     

At the same time, the heuristic study \cite{Gi} of $S^1$-equivariant 
geometry on the loop space $LX$ suggests that 
the generating functions $S=S_{1,\b}(0, Q, q)$ should satisfy certain
linear $q$-difference equations (instead of similar linear differential
equations of quantum cohomology theory). This expectation is
supported by the example of $X=\CC P^n$:
Y.-P. Lee \cite{Le} finds that the generating functions are solutions to 
the $q$-difference equation $D^{n+1} S=Q S$ (where $(D S) (Q):=S(Q)-S(qQ)$).  
 
In the case of the flag manifold $X$ the generating functions $S$ have 
been identified with the so called {\em Whittaker functions}
 --- common eigen-functions of commuting operators
of the $q$-difference Toda system. This result and its
conjectural generalization \cite{GL} 
to the flag manifolds $X=G/B$ of complex simple 
Lie algebras links quantum K-theory to representation theory and quantum
groups.  Originally this conjecture served as a motivation for
developing the basics of quantum K-theory.

\medskip

The author is thankful to Yuan-Pin Lee and Rahul
Pandharipande for useful discussions, and to 
the National Science Foundation for supporting this 
research.

\enddocument